\documentclass{m2an}
\usepackage{hyperref}
\begin{document}
%%-----------------------------
%%      the top matter
%%-----------------------------
\title{Nonlocal gradient operators with a nonspherical interaction neighborhood and their applications}\thanks{This 
	work is supported in part by the U.S.~NSF DMS-1719699, AFOSR MURI center for material failure prediction through peridynamics, and the ARO MURI Grant W911NF-15-1-0562.}
\author{Hwi Lee} \address{Department of Applied Physics and Applied Mathematics, Columbia University, New York, NY
	10027, USA.\\ \email{hl3001@columbia.edu; qd2125@columbia.edu}}
\author{Qiang Du} \sameaddress{1}
\date{...}
\begin{abstract} Nonlocal gradient operators are prototypical nonlocal differential operators that  are 
	very important in the studies of nonlocal models. One of the simplest variational settings for such studies is the nonlocal Dirichlet energies wherein the energy densities are quadratic in the nonlocal gradients. There have been earlier studies to illuminate the link between the coercivity of the Dirichlet energies and the interaction strengths of radially symmetric kernels that constitute nonlocal gradient operators in the form of integral operators. In this work we adopt a different perspective and focus on nonlocal gradient operators with a non-spherical interaction neighborhood. We show that the truncation of the  spherical interaction neighborhood  to a half sphere helps making nonlocal gradient operators well-defined and 
	the associated nonlocal Dirichlet energies coercive. These become possible, 
	unlike the case with full spherical neighborhoods,  without any extra assumption on the strengths of the  kernels
	near the origin.  We then present some applications of the nonlocal gradient operators with non-spherical interaction neighborhoods. These include nonlocal  linear models in mechanics such as nonlocal isotropic linear elasticity and nonlocal Stokes equations, and a nonlocal extension of the Helmholtz decomposition. \end{abstract}
%
%\begin{resume} ... \end{resume}
%
\maketitle
\nocite{*}
%%-----------------------------
%%      your text
%%-----------------------------
\section*{Introduction}
Nonlocal continuum models have become increasingly popular in many  scientific fields  \cite{AMRT10,auko09,BBM01,buades2010image,bv16bk,burago2014graph,BuLe11,coifman2006diffusion,du15fcaa,DeMa89,fkk03nonlocal,go08,KLS10,li2017quasi,pismen01,s00,tbl06,ToEn03,ts2016var,berto12gamma,zd10}. In particular nonlocal models given in terms of integral equations have received much attention as  alternatives, tools and bridges to classical local differential equation models  and discrete models; see \cite{du19} for more references and further discussions. Nonlocal models can be particularly effective in modeling singular physical phenomena. Peridyanmics, for example, was proposed by Silling \cite{s00} as a continuum theory to study  materials failure to which classical continuum theories of elasticity are not well suited. Central to the nonlocal models are nonlocal operators that are integral relaxations of the counterpart local differential operators, and this work is concerned with a class of representatives of such nonlocal operators, namely nonlocal gradient operators. 

Nonlocal gradient operators have been studied in various contexts ranging from rigorous theoretical studies on associated nonlocal function spaces to applications of the operators to nonlocal a posteriori stress analysis \cite{md16,dtty16}. The development of nonlocal vector calculus \cite{dm13} has also provided a systematic foundation to the study of nonlocal gradient operators. What motivates our work here are in two folds: the most relevant one is  the stability of nonlocal systems associated with nonlocal gradient operators as discussed in the context of correspondence theory of peridynamics \cite{dt17} and as utilized in the setting of fluid dynamics \cite{dtsph17}. Du and Tian \cite{dt17,dtsph17} have established the stability provided that nonlocal interaction kernels are radially symmetric and have suitably strong singularity of fractional type at the center of the nonlocal interaction neighborhood. We extend their analysis by proving that the use of nonspherical interaction neighborhood (or non-radially defined nonlocal interaction kernel) allows the assumptions on singular interaction  at the center to be removed when establishing the coercivity of the nonlocal Dirichlet energies. The idea of breaking the radial symmetry in the nonlocal interaction contributes to the second motivation of our work in connection to the study of such nonlocal operators for nonlocal convections \cite{tjd17} and nonlocal in time dynamic processes \cite{dyz17}. The former is important to preserve the upwinding feature of the physical transport process while the latter reflects the time-irreversibility.

Existing studies in the literature have demonstrated successful applications of nonlocal operators with non-radially symmetric kernels to nonlocal modeling. In \cite{dyz17}, the well-posedness and localization of nonlocal in time parabolic equation is established for a wide class of kernels when the support of the kernels for nonlocal-in-time derivative operators are truncated to yield one-sided backward in time nonlocal derivative operators. In study of compact embeddings of nonlocal function spaces of $L^p$ vector fields \cite{dmt18}, it is shown that the monotonicity assumption of the kernels for the nonlocal diffusion operators can be relaxed if they are non-radially symmetric and remain non-negative on a conic region. In \cite{fkv15}, the unique solvability is examined for nonlocal diffusion operators with kernels that may not have any symmetry at all. One of our main contributions here is to establish the coercivity result on the  correspondence  nonlocal Dirichlet energies provided that, for nonlocal gradient operators used,  one starts with radially symmetric kernels and breaks the symmetry by truncating the support of the kernels via multiplication by characteristic functions associated with half-spheres. 

The organization of the paper is as follows. 
In Section 2 we propose nonlocal nonsymmetric gradient operators and verify their consistency in the case of linear functions. We then introduce associated adjoint operators, namely nonlocal divergence operators, as well as corresponding nonlocal diffusion operators. We present the spectra of the nonlocal operators introduced thus far by imposing the periodic boundary conditions. This allows us to furnish, without technicalities of nonlocal boundary conditions, a straightforward proof of the main result of the paper that our nonlocal diffusion operators are positive definite \emph{uniformly} in two nonlocality parameters: one for the length scale and the other for the geometric specifications of nonlocal interaction neighborhoods. We conclude the Section 2 with discussion of our nonlocal diffusion operators in relation to other existing formulations of the operators, followed by a brief discourse on the issue of uncountably many choices of nonlocal interaction neighborhoods in our multi-dimensional formulation. In Section 3 we demonstrate the applications of our nonlocal gradient operators within the purview of linearized local continuum mechanics models. Specifically we revisit the unique solvability of the nonlocal Stokes equation \cite{dtsph17} and provide a nonlocal version of the classical Helmholtz decompositions. We also study a nonlocal model for isotropic linear elasticity. In Section 4 we give concluding remarks.

\section{Nonlocal nonsymmetric gradient operators}

A general form of the nonlocal gradient operator $\mathcal{G}_\delta$ studied in earlier works \cite{du19,dm13,dtty16,md16} is given by
\begin{equation} \label{ngg}
\mathcal{G}_\delta u(x) = 2\int_{\mathbb{R}^d}  \underline{\omega}_\delta(y-x) (y-x)\otimes \frac{u(y)-u(x)}{|y-x|} dy
\end{equation}
for a suitably defined function $u=u(x): 
{\mathbb{R}^d} \to \mathbb{R}^{1} \text{ or } \mathbb{R}^d$, and a nonlocal interaction kernel $\underline{\omega}_\delta=   \underline{\omega}_\delta(z)$. A common case for the latter is given by a radially symmetric kernel.  In this work, we define a nonlocal gradient operator $\mathcal{G}^{\vec{n}}_\delta$ acting on $u$ as  
\begin{equation} \label{ng}
\mathcal{G}^{\vec{n}}_\delta u(x) = 2\int_{\mathbb{R}^d} {\chi}_{\vec{n}}(y-x) w_\delta(|y-x|) (y-x)\otimes \frac{u(y)-u(x)}{|y-x|} dy
\end{equation}
where ${\chi}_{\vec{n}}(z)$ denotes the characteristic function of the half-space $\mathcal{H}_{\vec{n}} = \{\vec{z}\in \mathbb{R}^d:  \vec{z}\cdot \vec{n}\geq 0 \}$ parameterized by the unit vector $\vec{n}$.  Here,  the scalar-valued nonnegative function $w_\delta$ is a radially symmetric nonlocal kernel that measures the strength of the nonlocal interaction.   
We consider in particular a scaled kernel of the form $w_\delta(|x|) = \frac{1}{\delta^{d+1}}w(\frac{|x|}{\delta})$ where $w$ is nonnegative and compactly supported on the unit ball with a bounded first moment
\begin{equation} \label{mc}
\int_{\mathbb{R}^d} w(|x|)|x| dx= \int_{|x|\leq 1} w(|x|)|x| dx= d.
\end{equation}

Corresponding to \eqref{ngg}, we have $ \underline{\omega}_\delta(z)={\chi}_{\vec{n}}(z) w_\delta(|z|)$, which is no longer radially symmetric. The support of $w_\delta$ is given by the ball of radius $\delta$, which is called a nonlocal horizon or smoothing length depending on different contexts \cite{s00,dtsph17}.  Effectively, the nonlocal interaction neighborhood (the domain of integration) in \eqref{ng} is given by the support of $\underline{\omega}_\delta$ inside a half sphere defined by $\vec{n}$ with radius $\delta$.
The operators $\mathcal{G}_\delta^{\vec{n}}$ are $d$-dimensional versions of the following one dimensional one-sided nonlocal gradient operators  \cite{dtty16,dtsph17,dyz17}
\begin{equation} \label{1ng}
\mathcal{G}_\delta^{\pm}u(x) =  \pm 2 \int_{0}^\delta \omega_\delta(s)(u(x\pm s))-u(x))ds.
\end{equation}
Nonlocal derivative operators with a nonradial interaction kernel or a nonspherical neighborhood have also been used to study spatially inhomogeneous nonlocal convection in multidimensional cases \cite{tjd17}. However, systematic studies remain limited. We thus present some further analysis in this section.

\subsection{Consistency for linear functions}

Let us first note that the operators $\mathcal{G}_\delta^{\vec{n}}$ coincide with their local counterparts for linear functions.

\begin{lmm}
	\label{lem1}
	For every unit vector $\vec{n}$ in $\mathbb{R}^d$ and every affine function $u: \mathbb{R}^d \to \mathbb{R}^{\tilde{d}}$ where $\tilde{d} = 1 \text{ or } d$,
	$$
	\mathcal{G}_\delta^{\vec{n}} (u)(x) = \nabla u(x) \quad \forall x \in \mathbb{R}^d.
	$$	
\end{lmm}
\begin{proof}
	We assume without loss of generality $\tilde{d} = d$ and $x = 0$. Let us write $u(x) = Ax + b$ for some $d \times d$ real-valued matrix $A$ and $b \in \mathbb{R}^d$. If $R$ is a rotation matrix that aligns $\vec{e}_1 =(1,0,\dots,0) \in \mathbb{R}^d$ with $\vec{n}$, then we have
	\begin{align*}
	\mathcal{G}_\delta^{\vec{n}} (u)(0) &= 2 \int_{\mathcal{H}_{\vec{n}}}w_\delta(|y|) y \otimes \frac{Ay}{|y|}dy = 2 \int_{\mathcal{H}_{\vec{e}_1}}w_\delta(|z|) Rz \otimes \frac{ARz}{|z|}dz \\ &
	= R \left(2\int_{\mathcal{H}_{\vec{e}_1}}w_\delta(|z|) z \otimes \frac{z}{|z|}dz \right) R^TA^T = RI_dR^TA^T = A^T=\nabla u (0)
	\end{align*}
	where $T, I_d$ denote the transpose operator and the $d \times d$ identity matrix, respectively, and the moment condition (\ref{mc}) is used in the third to the last equality.
\end{proof}

\subsection{Adjoint operator}

In analogy with the local setting, the definition of the operator $\mathcal{G}_\delta^{\vec{n}}$ leads to the consideration of its adjoint  nonlocal divergence operator $\mathcal{D}^{\vec{n}}_\delta$. That is,  we define the operator $\mathcal{D}^{\vec{n}}$ by
\begin{equation} \label{anib}
\int_{\mathbb{R}^d} u(x)\cdot \mathcal{G}^{\vec{n}}_\delta v(x)dx = -\int_{\mathbb{R}^d}  \mathcal{D}_\delta^{\vec{n}}  u(x) v(x)dx
\end{equation}
for suitable functions $u$ and $v$ such that both sides of the equality make sense. In explicit terms, the operator $\mathcal{D}_\delta^{\vec{n}}$ takes the form
$$
\mathcal{D}_\delta^{\vec{n}} u(x)= 2 \int_{\mathbb{R}^d} {\chi}_{\vec{n}}(s) \frac{w^\delta(|s|)}{|s|}s\cdot ({u(x)-u(x-s)})ds.
$$

Note that in an analogy to the local diffusion operator $\Delta = \text{div grad}$, we may also define the associated nonlocal diffusion operator $\mathcal{L}_\delta^{\vec{n}} = \mathcal{D}^{\vec{n}}_\delta \circ \mathcal{G}^{\vec{n}}_\delta$ where $\circ$ denotes the composition. More discussions on
$\mathcal{L}_\delta^{\vec{n}}$ are given later.

\subsection{Representation in Fourier space}
For simplicity and definiteness, we assume from here on the two dimensional setting $d=2$ unless otherwise noted. Moreover we focus only on the action of the nonlocal operators on functions that are periodic on the domain $\Omega = (-\pi, \pi)^2$. We can then exploit the periodicity to examine the Fourier symbols of the nonlocal operators introduced thus far. In particular, for any periodic function $u(x)$ on $\Omega$, we let $\widehat{u}=\widehat{u}(\xi)$ denote its Fourier coefficients, hence the Fourier expansion
$$u(x)=\sum_{\xi\in\mathbb{Z}^2} \widehat{u}(\xi) e^{i\xi\cdot x} .$$
The loss of symmetry in the integration domain in the definitions of the nonlocal operators,
in contrast to the symmetric case studied in \cite{dtsph17},  is manifested in terms of the real parts of the Fourier symbols.

\begin{lmm} \label{l1}
	The Fourier symbols of the operators $\mathcal{G}_{\delta}^{\vec{n}},\mathcal{D}_{\delta}^{\vec{n}},\mathcal{L}_{\delta}^{\vec{n}}$ are given by
	\begin{align*}
	\widehat{\mathcal{G}_\delta^{\vec{n}} u} (\xi) &= \lambda_\delta^{\vec{n}}(\xi) (\widehat{u}(\xi))^T \\
	\widehat{\mathcal{D}_\delta^{\vec{n}} v} (\xi) &= -\overline{\lambda_\delta^{\vec{n}}(\xi)}^T \widehat{v}(\xi) = 
	%\textcolor{black}
	{{\lambda_\delta^{-\vec{n}}(\xi)}^T \widehat{v}(\xi)} \\
	\widehat{\mathcal{L}_\delta^{\vec{n}} w} (\xi) &= -|\lambda_\delta^{\vec{n}}(\xi)|^2 \widehat{w}(\xi)
	\end{align*}
	where $\xi \in \mathbb{Z}^2$ and
	\begin{align} \label{fng}
	\Re(\lambda_\delta^{\vec{n}}(\xi)) &=  2 \int_{\mathcal{H}_{\vec{n}}} \frac{w_\delta(|s|)s}{|s|}(\cos(\xi\cdot s)-1) ds \\ 
	\Im(\lambda_\delta^{\vec{n}}(\xi)) &= 2\int_{\mathcal{H}_{\vec{n}}} \frac{w_\delta(|s|)s}{|s|}\sin(\xi \cdot s ) ds.
	\end{align}
\end{lmm}

The above results are immediate from the definitions of the operators. It is natural to compare
$\lambda_{\delta}^{\vec{n}}(\xi)$ to the Fourier symbol of the local gradient operator given by
$ i  {\xi}$ to relate the imaginary part of $\lambda_{\delta}^{\vec{n}}(\xi)$ to %\textcolor{black}
{its local counterpart}. The former is shown to be a scalar multiple of the latter independently of ${\vec{n}}$ due to some  symmetry of the integrand on half-spheres.

\begin{lmm} \label{l2}
	For each ${\vec{n}}$, the Fourier symbol $\lambda_\delta^{\vec{n}}(\xi)$ in (\ref{fng}) can be expressed as
	\begin{equation}
	\lambda_\delta^{\vec{n}}(\xi) = i \Lambda_\delta(|\xi|)\frac{\xi}{|\xi|} +  \Re(\lambda_\delta^{\vec{n}}(\xi))
	\end{equation}
	where
	\begin{equation}
	\Lambda_\delta(|\xi|) = 4 \int_{s_1 \geq 0, s_2 \geq 0}  \frac{w_\delta(|s|)s_1}{|s|}\sin(|\xi|s_1)ds.
	\end{equation} 
	%	\textcolor{black}
	{On the other hand, $\Re(\lambda_\delta^{\vec{n}}(\xi))$ is a scalar multiple of $\frac{\xi}{|\xi|}$ if and only if $\vec{n}$ is a scalar multiple of $\xi$.}
\end{lmm} 

\begin{proof}
	We observe from Lemma \ref{l1}  that
	$$
	\Im(\lambda_\delta^{\vec{n}}(\xi)) = 2 \int_{\mathcal{H}_{\vec{n}}} \frac{w_\delta(|s|)s}{|s|}\sin(\xi \cdot s ) ds = \int_{B_\delta(0)} \frac{w_\delta(|s|)s}{|s|}\sin(\xi \cdot s ) ds
	$$
	where the equality is due to the odd symmetry of $s \sin(\xi\cdot s)$. Then the first claim follows from Lemma 3 in \cite{dtsph17}.
	Next we consider the real part of $\lambda_\delta^{\vec{n}}(\xi)$. If we let $\xi^{\perp}$ denote a vector orthogonal to $\xi$, then it follows that
	$$
	|\xi^{\perp} \cdot \Re(\lambda_\delta^{\vec{n}}(\xi)| = 2 \int_{\mathcal{I}_{\vec{n},\xi} } \frac{w_\delta(|s|)|\xi^{\perp} \cdot s|}{|s|}(1-\cos(\xi \cdot s))ds
	$$
	where $\mathcal{I}_{\vec{n},\xi} = \{s \in \mathcal{H}_{\vec{n}}: s - 2 (s \cdot \xi) \xi  \in \mathcal{H}_{\vec{n}}\}$, thus the second claim holds since $|\mathcal{J}_{\vec{n},\xi}| = 0$ precisely when $\vec{n}$ is a scalar multiple of $\xi$.
\end{proof}

\subsection{Spectral estimates}

We now present a key theorem concerning the spectral property 
of the nonlocal gradient operator $\mathcal{G}_\delta^{\vec{n}}$. 
The theorem implies, in particular the strong coercivity, uniformly in $\delta$ and $\vec{n}$, of the Dirichlet energies given by {$\| \mathcal{G}_\delta^{\vec{n}}(\cdot) \|_2^2$ 
	with respect to the norm $\| \cdot \|_{2} + \| \mathcal{G}_\delta^{\vec{n}}(\cdot) \|_2.$}

\begin{thrm} \label{t1}			
	{There exists a positive constant $C$ independent of $\vec{n}$, $\xi$ and $\delta$ (as $\delta\to 0$) such that 
		$$C \leq |\lambda_\delta^{\vec{n}}(\xi)|\leq \sqrt{2} d  |\xi| \quad \forall \, \xi \in \mathbb{Z}^2/\{0\}.$$}
\end{thrm}

\begin{proof}
	Let $k = \delta |\xi|$, we show the bound on $|\lambda_\delta^{\vec{n}}(\xi)|$
	using two separate estimates on the imaginary and real parts of
	$\lambda_\delta^{\vec{n}}(\xi)$ respectively depending on $k< 1$ or $k\geq 1$. 	
	
	For $k < 1$, we estimate the imaginary part.
	Using Lemma \ref{l2}, we have
	\begin{align*}
	\Lambda_\delta(\xi)  = 4 \int_{s_1 \geq 0, s_2 \geq 0}\frac{w_\delta(|s|)s_1}{|s|}\sin(|\xi|s_1)ds 
	& = 4 \int_{r=0}^{\delta}\int_{\theta = 0}^{\frac{\pi}{2}} w_\delta(r) r\cos(\theta) \sin(|\xi|r \cos(\theta))drd\theta \\
	& = \frac{4}{\delta} \int_{r=0}^{1}\int_{\theta=0}^{\frac{\pi}{2}} w(r)r\cos(\theta)\sin(|\xi|\delta r \cos(\theta))dr d\theta\,.
	\end{align*}
	By the inequality
	$
	\sin(x) \geq x - {x^3}/{6}$  for  $0 \leq x \leq 1$,
	we obtain
	\begin{align*}
	\Lambda_\delta(\xi) &\geq \frac{4k}{\delta}\int_{r=0}^1\int_{\theta = 0}^{\frac{\pi}{2}}w(r)r^2\cos^2(\theta)drd\theta - \frac{4k^3}{6\delta}\int_{r=0}^{1}\int_{\theta = 0}^{\frac{\pi}{2}}w(r)r^4\cos^4(\theta)drd\theta \\
	&\geq \frac{C k}{\delta} = C |\xi|,
	\end{align*}
	for a constant $C>0$, independent of $\xi$, $\delta$ and $\vec{n}$.
	
	Next, for $1 \geq k$, we consider $\Re (\lambda_\delta^{\vec{n}}(\xi))$. Note that
	\begin{align*}
	|\Re (\lambda_\delta^{\vec{n}}(\xi))| \geq |\vec{n} \cdot \Re (\lambda_\delta^{\vec{n}}(\xi))| &= 2 \int_{\mathcal{H}_{\vec{n}}} \frac{w_\delta(|s|)}{|s|}(n_1 s_1 + n_2 s_2)(1-\cos(\xi \cdot s))ds \\
	&\geq 2\cos\left(\frac{\pi}{4}\right) \int_{H_{\vec{n},\frac{\pi}{4}}} w_\delta(|s|) (1-\cos(\xi \cdot s))ds
	\end{align*}
	where ${\mathcal{H}_{\vec{n},\frac{\pi}{4}}}$ denotes the set of those points in $\mathcal{H}_{\vec{n}}$ that have 
	angles with $\vec{n}$ between {$-\pi/4$} 
	and $\pi/4$. In terms of polar coordinates, we can then write
	\begin{align*}
	|\Re (\lambda_\delta^{\vec{n}}(\xi))|  & \geq C \int_{r = 0}^{\delta} \int_{\theta_1}^{\theta_1 + \frac{\pi}{2}} w_\delta(r) (1-\cos(|\xi|r \cos(\psi_{\theta})))rdrd\theta\\
	& = \frac{C}{\delta} \underbrace{\int_{r = 0}^{1} \int_{\theta_1}^{\theta_1 + \frac{\pi}{2}} w(r) (1-\cos(kr |\cos(\psi_{\theta})|))rdrd\theta}_{J_\delta(\xi)}
	\end{align*}
	for some $\theta_1$ depending on $\vec{n}$. Here  $0 \leq \psi_\theta \leq \pi$ denotes the angle between the vector $\xi$ and the vector with the polar coordinates $(r,\theta)$. 
	
	We now introduce a possible cut-off  of $w$ at the origin to get an absolutely integrable kernel $\phi$,
	that is, $\phi$ is a radial function such that
	$$0 \leq \phi(|x|) \leq w(|x|) \text { and } 0 < I := \int_{\mathbb{R}^2} \phi(|x|)dx < \infty. $$

	We then discuss the different cases separately.  First, let us consider
	the case with $1 \leq k \leq \lambda$ where $\lambda$ is to be specified.

	If we let 
	$$
	A_\xi =  \left(\theta_1,\theta_1 + \frac{\pi}{2}\right) - \left\{\theta \in \left(0,2 \pi \right): \left|\psi_\theta - \frac{\pi}{2} \right| \leq \frac{\pi}{8} \right\},
	$$
	then 
	$$ \frac{\cos\left(\frac{3\pi}{8}\right)}{2 \lambda}\leq kr |\cos(\psi_\theta)| \leq 1 \quad  \text{ for }  (r,\theta) \in \left(\frac{1}{2\lambda},\frac{1}{\lambda}\right) \times A_\xi$$ so that
	\begin{align*}
	J_\delta(\xi) &\geq \frac{1}{\delta} \int_{r = \frac{1}{2\lambda}}^{\frac{1}{\lambda}} \int_{A_\xi} w(r)r(1-\cos(k r |\cos(\psi_\theta)|))dr d\theta\\
	& \geq \frac{1}{\delta}\left(1-\cos\left(\frac{\cos\left(\frac{3\pi}{8}\right)}{2 \lambda}\right)\right) \int_{r = \frac{1}{2\lambda}}^{\frac{1}{\lambda}} \int_{A_\xi} w(r)rdr d\theta \geq \frac{C}{\delta}
	\end{align*}
	where the last inequality is due to the non-degeneracy of $|A_\xi| \geq \frac{\pi}{4}$ uniformly in $\xi$. 
	
	Next, we consider the case $\lambda < k$.
	Then with the same $A_\xi$ defined in the case (a) 
	\begin{align*}
	J_\delta(\xi) &\geq \frac{1}{\delta}\int_{A_\xi} \int_{r=0}^{1} w(r)r \left(1-\cos\left(kr|\cos(\psi_\theta)|\right)\right) dr d\theta \\
	&\geq \frac{1}{\delta} \int_{A_\xi} \left(\frac{I}{2\pi} - \int_{r=0}^1 \phi(r)r \cos(kr |\cos(\psi_\theta)|)dr \right) d\theta.
	\end{align*}
	{By the Riemann-Lebesgue lemma, there exists a constant $c > 0$ such that for $j > c$
		$$
		\left|\int_{r=0}^1 \phi(r)r \cos\left(jr \right)dr\right| < \frac{I}{4 \pi}.
		$$
		Then, since $\cos(\frac{3\pi}{8}) \leq |\cos(\psi_\theta)|$ for $\theta \in A_\xi$ we set $\lambda = c/{\cos(\frac{3\pi}{8})}$ to obtain that for some $\tilde{C}$,} 
	$$
	J_\delta(\xi) \geq \frac{\tilde{C}}{\delta}.
	$$		
	In summary, we have $|\lambda_\delta^{\vec{n}}(\xi)| \geq \min\{{C_1},\frac{C_2}{\delta} \}$ for positive constants $C_1$ and $C_2$.
	
	In order to prove the uniform upper bound on $|\lambda_\delta^{\vec{n}}(\xi)|$ we observe
	{			\begin{equation*} 
		|\Re(\lambda_\delta(\xi))| =
		\left| 2 \int_{\mathcal{H}_{\vec{n}}} \frac{w_\delta(|s|)s}{|s|}(\cos(\xi\cdot s)-1) ds \right| \leq
		2 \int_{\mathcal{H}_{\vec{n}} } 
		{w_\delta(|s|)}   |\xi \cdot s | ds \leq d |\xi|
		\end{equation*}	
		and
		\begin{equation*}
		|\Im(\lambda_\delta(\xi))| = \left|  \int_{B_\delta(0)}  \frac{w_\delta(|s|)s}{|s|}\sin(\xi \cdot s)ds\right| 
		\leq  \int_{B_\delta(0)}  w_\delta(|s|)|\xi| |s|ds 
		= d  |\xi|
		\end{equation*}
		following the derivation in the  proof of Lemma \ref{l2}.}
\end{proof}

\subsection{Nonlocal  diffusion and nonlocal gradient}
With the assurance of the coercivity we next address the question of whether there is a connection between the nonlocal diffusion operator based on the non-symmetric formulations of nonlocal divergence and gradient and other existing formulations of nonlocal diffusion (Laplacian) operator. We begin by recalling from \cite{dt17} that symmetric gradient operator is closely related to a bond-based, in the language of peridynamics, nonlocal diffusion operator. We show that the similar conclusion can be drawn in our formulation with some choices of kernels. 

To illustrate this point it is sufficient to consider the one dimensional setting wherein the nonlocal Dirichlet integrals are 
$$
\mathcal{E}_\delta^{\pm}(u) = \int_\Omega |\mathcal{G}^{\pm}_\delta(u)(x)|^2 dx.
$$
Since we know from Lemma \ref{l2} that $|\lambda_\delta^{\vec{n}}(\xi)| = |\lambda_\delta^{-\vec{n}}(\xi)|$, it is sufficient to consider only $\mathcal{E}^{+}$ for which we have the following. 

\begin{lmm} \label{cl}
	If $w_\delta(|x|)$ is integrable, then
	$$
	{\mathcal{E}}^{+}_\delta(u) =  2 \int_\Omega \int_{0}^{\delta} {\rho}_{\delta}(a)\left|\frac{u(x+a)-u(x)}{a}\right|^2 dadx
	$$
	where ${\rho}_\delta(a)={\rho}_\delta(|a|)$ is a radial (even) function given by 
	\begin{align}
	\label{rhodelta}
	{\rho}_\delta(a) &= 2a^2 \int_{0}^{\delta} w_\delta(b)\left(w_\delta(a)-w_\delta\left(a+b\right)\right)db
	= - a^2  (\mathcal{G}^{+}_\delta w_\delta) (a) ,\quad \forall a\in (0, \delta)
	\end{align}
	supported on $(-\delta,\delta)$ and satisfies $\|\rho_\delta\|_{L^1} = 1$. Moreover if $\omega$ is non-increasing in $(0,1)$, then $\rho_\delta$ is non-negative.
\end{lmm}

\begin{proof}
	Following the argument in \cite{dt17}, we write $D_\delta = [-\delta,0]^2 \bigcup [0,\delta]^2$ so that
	\begin{align*} 
	%\textcolor{black}
	{\mathcal{E}_\delta^{+}}(u) &= 2 \underbrace{\int_\Omega \int_{D_\delta} w_\delta(|s|)w_{\delta}(|t|) (u(x+t)-u(x))^2 dsdtdx}_{I_1} + \dots \nonumber \\ & \qquad \quad \underbrace{-\int_\Omega \int_{D_\delta} w_\delta(|s|)w_{\delta}(|t|) ((u(x+s)-u(x+t))^2 dsdtdx}_{I_2}.
	\end{align*}
	Let us first observe
	\begin{equation*}
	I_1 = \int_\Omega \int_{-\delta}^{\delta} k_{\delta}(|a|) \left|\frac{u(x+a)-u(x)}{a}\right|^2 dadx
	\end{equation*}
	where 
	$$
	k_\delta(|a|) = 2w_{\delta}(|a|)|a|^2\int^{\delta}_{0} w_\delta(|b|)db.
	$$
	Next we consider $I_2$ and rewrite it as
	\begin{align*}
	I_2 &= -\int_\Omega \int_{D_\delta} w_\delta(|s|)w_{\delta}(|t|) ((u(x+s-t)-u(x))^2 dsdtdy \\
	&= -\int_\Omega \int_{\hat{D}_\delta}\frac{a^2}{2} w_\delta\left(\left|\frac{a+b}{2}\right|\right)w_{\delta}\left(\left|\frac{b-a}{2}\right|\right) \frac{((u(x+a)-u(x))^2}{a^2} dadbdx
	\end{align*}
	where we have used the periodicity of $u$ in the first equality and the change of variables $a = s-t$ and $b = s+t$ in the second with the corresponding change of integration domain from $D_\delta$ to $\hat D_\delta$. That is we have
	$$
	I_2 = \int_\Omega \int_{-\delta}^{\delta} h_\delta(|a|) \frac{(u(y+a)-u(y))^2}{a^2} dadx
	$$
	where
	\begin{align*}
	h_\delta(|a|) &= -\frac{a^2}{2} \left(\int_{|a|}^{-|a|+2\delta} + \int_{|a|-2\delta}^{-|a|} \right)  w_\delta(|\frac{a+b}{2}|)w_{\delta}(|\frac{b-a}{2}|) db = -{a^2} \int_{|a|}^{-|a|+2\delta}  w_\delta(|\frac{a+b}{2}|)w_{\delta}(|\frac{b-a}{2}|) db. 
	\end{align*}
	We note the almost everywhere finiteness of $h_\delta(|a|)$ since
	$$
	\int_{-a}^{a} h_\delta(|a|) da = 1 - \left(\int_{-\delta}^{\delta} s^2w_\delta(|s|)ds\right)\left(\int_{-\delta}^{\delta} w_\delta(|t|)dt\right)
	$$
	which can be easily verified. Then the desired claim follows from setting ${\rho}_\delta = k_\delta + h_\delta$ and observing that for $0< a < \delta$,
	\begin{align*}
	h_\delta(|a|) &= -{2a^2} \int_{0}^{\delta-a}  w_\delta(|z+a|)w_{\delta}(|z|) dz = -{2a^2} \int_{0}^{\delta}  w_\delta(|z+a|)w_{\delta}(|z|) dz
	\end{align*}
	where the last equality holds since $w_\delta$ is supported on $(-\delta,\delta)$, and analogously for $ -\delta < a < 0$.
	
	Finally it is clear from (\ref{rhodelta}) that $\rho_\delta$ is non-negative for non-increasing $\omega_\delta$.
\end{proof}

Naturally, one may be interested in extending the result of the above lemma to more general
kernels $w$.  Let us consider first a special example
$$
w(|x|)=
\begin{cases}
\frac{C_\beta}{|x|^\beta}, \quad |x| \leq 1,\\
0, \quad \text{ otherwise}
\end{cases}
$$
for $1\leq \beta < 2$, and $C_\beta > 0$ is chosen to satisfy the moment condition (\ref{mc}). Then
the result of Lemma \ref{cl} also holds. Indeed,
let us fix $0 < a < \delta$ without loss of generality. We then have 
\begin{align*}
\frac{{\rho}_\delta(a)}{2a^2 C_\beta^2} &=  \int_0^{\delta-a} \frac{1}{z^{\beta}}\left(\frac{1}{a^{\beta}} - \frac{1}{(a+z)^\beta} \right)dz + \int_{\delta-a}^{\delta}\frac{1}{z^{\beta}}\frac{1}{a^{\beta}}dz \leq \int_0^{\delta-a} \frac{1}{z^{\beta}}\left(\frac{\beta  \delta^{\beta-1} z}{a^\beta(a+z)^\beta} \right)dz + \int_{\delta-a}^{\delta}\frac{1}{z^{\beta}}\frac{1}{a^{\beta}}dz < \infty.
\end{align*}
Moreover if we define 
$$
\rho_\delta^{\epsilon}(a) = {\chi}_{(\epsilon,\infty)}(|a|)2a^2 \int_{\epsilon}^{\delta} w_\delta(b)\left(w_\delta(a)-w_\delta\left(a+\frac{a}{|a|}b\right)\right)db
$$ which are nonnegative, monotonically increasing in $\epsilon$, pointwise convergent approximations of $\rho_\delta$ as $\epsilon \to 0$, then a direct calculation shows  
$$\|\rho_\delta\|_{L^1(\mathbb{R})} = \lim_{\epsilon \to 0} \int_{-\delta}^{\delta}\rho_\delta^{\epsilon}(a)da  = 1$$ where the first equality is due to the Monotone Convergence Theorem.

We can extend the definition of $\rho_\delta$ to a more broader class of non-integrable kernels $w_\delta$ that include the above fractional ones as a special case. To this end let us assume that $w_\delta$ is non-increasing and consider
\begin{align}
\label{rhodelta1}
w^\epsilon_\delta(|x|)=
\begin{cases}
w_\delta(|x|), &\quad |x| > \epsilon\\
\inf_{|y| \leq \epsilon }w_\delta(|y|), & \text{ otherwise}
\end{cases}.
\end{align}
Note that  the modified nonnegative and radial kernel
$w^\epsilon_\delta$  is  integrable 
and also  nonincreasing in $(0, \delta)$. Then 
\begin{align}
\label{rhodelta2}
{\rho}^{\epsilon}_\delta(a) &= 2a^2 \int_{0}^{\delta} w^{\epsilon}_\delta(b)\left(w^{\epsilon}_\delta(a)-w_\delta^{\epsilon}\left(a+ \frac{a}{|a|}b\right)\right)db
\end{align}
is non-negative, monotonically increasing  in $\epsilon$  and satisfies $\lim_{\epsilon \to 0} \|\tilde{\rho}_\delta^{\epsilon}\|_{L^1(\mathbb{R})} = 1$. Hence we may define
\begin{align}
\label{rhodelta3}
\rho_\delta(a) = \lim_{\epsilon \to 0} {\rho}_\delta^{\epsilon}(a)
\end{align}
which satisfies $\|\rho_\delta\|_{L^1(\mathbb{R})} = 1$ due to the Monotone Convergence Theorem. We thus get
the following more general result.

\begin{lmm} \label{cl2}
	If $w_\delta(|x|)$ is non-increasing with bounded first moment, then
	$$
	{\mathcal{E}}^{+}_\delta(u) =  2 \int_\Omega \int_{0}^{\delta} {\rho}_{\delta}(a)\left|\frac{u(x+a)-u(x)}{a}\right|^2 dadx
	$$
	where ${\rho}_\delta(a)={\rho}_\delta(|a|)$ is a radial (even) function defined by \eqref{rhodelta1}-\eqref{rhodelta2}-\eqref{rhodelta3}.
\end{lmm}

We remark that the non-negativity of $\rho_\delta$ for non-increasing $w_\delta$ is in a clear contrast with the case of symmetric nonlocal gradient operators wherein the corresponding kernel is always sign changing. On the other hand it is easy to see that $\rho_\delta$ is not always non-negative and may change sign as in the case of
$$
{\rho}_\delta(x) = \pi x^2 \sin(\pi |x|) + \frac{\pi^2 x^2}{4}\left((|x|-1)\cos(\pi x) - \frac{\sin(\pi|x|)}{\pi}\right)
$$
which amounts to $w_\delta(x) = \frac{\pi}{2}\sin(\pi |x|)$ with $\delta = 1$.

\subsection{Nonlocal diffusion operator, a revisit}

As a corollary of the discussion in the preceding subsection, in one space dimension, 
we have the equivalence of the corresponding nonlocal diffusion operator  based on the
nonsymmetric gradient and divergence
$$
\mathcal{L}^{+}_\delta u(x) = (\mathcal{D}^{+}_\delta \circ \mathcal{G}^{+}_\delta )u(x) = 4 \int_0^{\delta} \int_{0}^{\delta} w_\delta(|s|)w_\delta(|t|)(u(x+s)-u(x)-u(x+s+t) +u(x+t))dsdt
$$
and the conventional bond based nonlocal diffusion operator 
$$
\mathcal{L}_{\delta}u(x) = 2\int_{-\delta}^{\delta} k_\delta(|a|) (u(x+a)-u(x)) da 
$$
upon setting $k_\delta(|a|) = \frac{{\rho}_\delta(|a|)}{a^2}$ for kernels discussed in the previous subsection.

In the case of an integrable $w_\delta$, we can further relate the operator $\mathcal{L}^{+}_\delta$  to  another nonlocal diffusion operator in a similar form, namely a doubly nonlocal Laplace operator $\mathcal{L}_{\delta,\epsilon}^{double}$ proposed in \cite{rw18}
$$
\mathcal{L}_{\delta,\epsilon}^{double}u(x) = \int_{-\delta}^{\delta} \int_{-\epsilon}^{\epsilon} \gamma_\delta(y)\eta_\epsilon(r)(u(x+y+r)-u(x)- u(x+r) +u(x+y))dydr.$$
Here $\gamma_\delta$ and $\eta_\epsilon$ are assumed to be radial, non-negative and compactly supported on $(-\delta,\delta)$ and $(-\epsilon,\epsilon)$, respectively. Further assumptions on the moments of the kernels are made, namely the normalized second moment of $\gamma_\delta$ and integrability of $\eta_\epsilon$ with unit mass. For clarity of comparison, we first observe 
$$
\mathcal{L}_{\delta,\epsilon}^{double}u(x) = (\mathcal{L}_{\delta} \circ \mathcal{A}_\epsilon)u(x) = (\mathcal{A}_\epsilon \circ \mathcal{L}_\delta)u(x)
$$ 
if we let $k_\delta = \gamma_\delta$ in $\mathcal{L}_\delta$ and define the averaging operator $\mathcal{A}_\epsilon$ by
$$
\mathcal{A}_\epsilon u(x) = \frac{1}{2}\int_{-\epsilon}^{\epsilon} \eta_\epsilon(|z|)(u(x+z)+u(x)) dz
= \frac{1}{2} u(x)  +  \frac{1}{2} \int_{-\epsilon}^{\epsilon} \eta_\epsilon(|z|)u(x+z)dz
.
$$

Evidently the operator $\mathcal{A}_\epsilon$ provides a simple averaging, and does not fundamentally alter the spectral properties of  $\mathcal{L}_{\delta,\epsilon}^{double}$. Thus, while involving an extra kernel, it does not change
the modeling capability overall.

\subsection{Nonlocal gradient with orientation dependence}

We recall that our formulation of the operators $\mathcal{G}_\delta^{\vec{n}}$ is motivated by a geometric consideration of generalizing the one dimensional half-spaces $(0,\infty)$ and $(-\infty,0)$ to their higher dimensional analogues. The $\vec{n}$ itself belongs to the non-trivial ambient space $S^{d-1} = \{\vec{n} \in \mathbb{R}^d: \|\vec{n}\|_2 = 1 \}$ for $d \geq 2$.  A specific choice of $\vec{n}$ leads to orientation dependence. 
It is a legitimate question to ask if such dependence is necessary (while maintaining a coercive Dirichret integral).
One can also ask how to pick
$\vec{n}$, in case that it is needed, in practice.  Since the
coercivity result of Theorem \ref{t1} is true for all $\vec{n}$, one possible approach is
to eliminate the dependence on $\vec{n}$ by defining a new energy functional in terms of the average over $\vec{n}\in S^{d-1}$. Indeed, if we take for concreteness a scalar valued $u$ in two dimensions, we may consider the averaged nonlocal Dirichlet integral
$$
{\mathcal{E}^{avg}_\delta}(u) = \frac{1}{2\pi} \int_{S^{1}} \int_{\Omega} | \mathcal{G}_\delta^{\vec{n}}u(x) |^2 dx dS.
$$ The coercivity of ${\mathcal{E}}_{avg}$ is immediate from that of each Dirichlet integral associated with $\mathcal{G}_{\delta}^{\vec{n}}$. 
One can view ${\mathcal{E}}_{avg}$ as a stabilized symmetric nonlocal Dirichlet integral $\mathcal{E}^{sym}_\delta(u) =  \int_{\Omega} | \frac{1}{2}({\mathcal{G}_\delta^{\vec{n}}+\mathcal{G}_\delta^{\vec{n}})u(x)}|^2 dx
$ since
$$
{\mathcal{E}^{avg}_\delta}(u) =  {\mathcal{E}^{sym}_\delta}(u) + \sum_{\xi \in  \mathbb{Z}^2/\{0\}} 	\frac{2}{\pi} \int_{{B_\delta}(0)}\int_{{B_\delta}(0)} \frac{w_\delta(|a|)w_\delta(|b|)a\cdot b}{|a||b|}\arcsin\left(\frac{a \cdot b}{|a||b|}\right) (\cos(\xi\cdot a)-1)(\cos(\xi\cdot b)-1)dadb. 
$$
Unfortunately, we are unable to express ${\mathcal{E}^{avg}_\delta}(u)$ as a Dirichlet integral of a nonlocal gradient operator. However, based on the above calculation and 
using a crude estimate ${\pi}x^2 \geq x \arcsin(x)$ we may alternatively consider a simpler looking, yet still coercive
$$
{\mathcal{E}}^{\star, \vec{k}}_\delta(u) = \int_{\Omega} | \mathcal{G}_\delta^{\vec{k}}u(x) |^2 dx
$$
where for any vector $\vec{k}\in \mathbb{R}^d$, $\mathcal{G}_\delta^{\star,\vec{k}}u(x)$ is a modified nonlocal gradient operator given by
\begin{align}
\label{newgra}
\mathcal{G}_\delta^{{\star},\vec{k}}u(x) = \int_{\mathbb{R}^d}w_\delta(|y-x|) \frac{y-x}{|y-x|} {(u(y)-u(x))}dy + \left(\int_{\mathbb{R}^d}w_\delta(|y-x|) {(u(y)-u(x))} dy \right)\vec{k}.
\end{align}
We remark that the adjoint operator of $\mathcal{G}_\delta^{{\star},\vec{k}}$ can be defined in a similar fashion as in subsection 2.2. This may prove to be useful in nonlocal modeling. 
Note that
$\mathcal{G}_\delta^{{\star},\vec{k}}$ can be seen as a special form of more general nonlocal gradient operators studied in 
\cite{md16}. It is also related to the more standard nonlocal gradient operator with a spherical interaction neighborhood corresponding to the form of  $\mathcal{G}_\delta^{{\star},\vec{k}}$  with $\vec{k}=\vec{0}$ (the zero vector). According to \cite{dt17,dtsph17}, the coercivity of the Dirichret integral corresponding
to this case, i.e., $\vec{k}=\vec{0}$,  depends on the choices of the kernel $w_\delta$.   
For nonzero  $\vec{k}$, we can have the coercivity of the Dirichret integral for
$\mathcal{G}_\delta^{{\star},\vec{k}}$ 
due to the second term 
in \eqref{newgra}.
This orientation dependent term is $O(\delta)$,  due to the moment condition in (\ref{mc}),
similar in spirit to how stability is attained in our recent work on the deterministic particle methods \cite{ld19}.

\section{Applications of the nonlocal gradient operators}
We now illustrate how the coercivity of the nonlocal Dirichlet energies can be utilized in several applications.
These include the application of our nonlocal operators as building blocks for an alternative formulation of the nonlocal Stokes equation studied in \cite{dtsph17}. Another application is to establish a nonlocal version of the Helmholtz decomposition. In addition,  the operators are used to construct well-defined models of nonlocal isotropic linear elasticity
that converge to the classical counterpart for any Poisson's ratio.

We use the nonsymmetric operator $\mathcal{G}_\delta^{\vec{n}}$ with a unit vector $\vec{n}$ defined in \eqref{ng}
for illustration, though similar discussions can be made for $\mathcal{G}_\delta^{\star,\vec{k}}$ with a constant nonzero
vector $\vec{k}\neq \vec{0}$ defined in \eqref{newgra}.
To avoid further technical complications, 
we assume from here on the kernels adopted in our nonlocal operators are positive almost everywhere. 

\subsection{Nonlocal Stokes Equation}
We {first} consider the {steady} nonlocal stokes equation 
\begin{eqnarray} \label{ns}
-\mathcal{L}^{\vec{n}}_\delta u_\delta^{\vec{n}} + \mathcal{G}_\delta^{\vec{n}} p^{\vec{n}}_\delta &=& f \nonumber \text{ in } \Omega  \\
-\mathcal{D}^{\vec{n}}_\delta u^{\vec{n}}_\delta &=& 0 \text{ in } \Omega
\end{eqnarray}
where $u^{\vec{n}}_\delta,p^{\vec{n}}_\delta, f$ are periodic functions on $\Omega$ and assumed to have zero means for compatibility. One of the motivation for studying nonlocal Stokes model is to better understand
methods like the smoothed particle hydronamics (SPH) \cite{GiMo:77,lucy:77}, see \cite{dtsph17} for more
references and discussions.

The nonlocal Stokes equation is obtained by applications of the nonlocal operators $\mathcal{L}_\delta^{\vec{n}},\mathcal{G}_\delta^{\vec{n}},\mathcal{D}_\delta^{\vec{n}}$ to their local counterparts in the classical Stokes equation
\begin{eqnarray} \label{ls}
-\Delta u + \nabla p &=& f, \nonumber \text{ in } \Omega  \\
-\nabla \cdot u &=& 0, \text{ in } \Omega.
\end{eqnarray}
In the Fourier space the system (\ref{ns}) can be written as  
\begin{equation} \label{fns}
A_\delta^{\vec{n}}(\xi)           
\begin{bmatrix} 
\widehat{u^{\vec{n}}_{\delta}}(\xi) \\  \widehat{p^{\vec{n}}_{\delta}}(\xi)
\end{bmatrix} = \begin{bmatrix} 
\widehat{f}(\xi) \\  0
\end{bmatrix}
\end{equation}
where  
\[
A_\delta^{\vec{n}} =
\begin{bmatrix}
|\lambda_\delta^{\vec{n}}(\xi)|^2 I_2 &  \lambda_\delta^{\vec{n}}(\xi) \\
\overline{\lambda_\delta^{\vec{n}}(\xi)}^T & 0
\end{bmatrix}.
\] 
The non-degeneracy of $\lambda_{\delta}^{\vec{n}}$ assured by Theorem \ref{t1} yields the well-posedness of (\ref{fns}).
\begin{lmm} \label{wns}
	The system (\ref{fns}) has a unique solution given by
	\begin{equation} \label{sfns}
	\begin{bmatrix}
	\widehat{u^{\vec{n}}_{\delta}}(\xi) \\  \widehat{p^{\vec{n}}_{\delta}}(\xi)
	\end{bmatrix} = (A_\delta^{\vec{n}})^{-1}\begin{bmatrix} 
	\widehat{f}(\xi) \\  0
	\end{bmatrix}
	\end{equation}
	where
	\[
	(A_\delta^{\vec{n}})^{-1} = 
	\begin{bmatrix} 
	\frac{1}{|\lambda_{\delta}^{\vec{n}}(\xi)|^2} (I - \frac{\lambda_{\delta}^{\vec{n}}(\xi)\overline{{\lambda_{\delta}^{\vec{n}}(\xi)}}^T  }{|\lambda_{\delta}^{\vec{n}}(\xi)|^2})   &      \frac{ {\lambda_{\delta}^{\vec{n}}(\xi)}}{|\lambda_{\delta}^{\vec{n}}(\xi)|^2}\\
	\frac{\overline{\lambda_{\delta}^{\vec{n}}(\xi)}^T}{|\lambda_{\delta}^{\vec{n}}(\xi)|^2}  & -1
	\end{bmatrix}.
	\] 
	Moreover, there exists a constant $C>0$ independent of $f$ and $\delta$, as $\delta \to 0$, such that
	\begin{equation}
	\|u^{\vec{n}}_\delta\|_{S_\delta^{\vec{n}}(\Omega)} + \|p^{\vec{n}}_\delta\|_{L^2(\Omega)} \leq  C \|f\|_{(S_\delta^{\vec{n}}(\Omega))^{\star}} 			
	\end{equation}
	where $S_\delta^{\vec{n}}(\Omega)$ is the energy space equipped the norm $(\sum |\lambda_{\delta}^{\vec{n}}(\xi)|^2|\widehat{u}(\xi)|^2)^{1/2}$, and  $(S_\delta^{\vec{n}}(\Omega))^{\star}$ is its dual space with respect to the standard $L^2$ duality pairing.  
\end{lmm}

Let us examine the limit of the nonlocal solutions $(u^{\vec{n}}_\delta,p^{\vec{n}}_\delta)$ as $\delta \to 0$. As can be expected by Lemma \ref{l2}, it is worthy noting that in the case studied in this work, the nonlocal velocity is only approximately local divergence free, which is in contrast with the nonlocal Stokes equation in \cite{dtsph17} that gives equivalent  local and nonlocal
divergence free vector fields.

\begin{prpstn} \label{p1}
	Let $(u,p)$ and $(u^{\vec{n}}_\delta,p^{\vec{n}}_\delta )$ denote the solutions of (\ref{ns}) and (\ref{ls}) respectively. Then there exists a constant $C>0$ independent of $f$ and $\delta$, as $\delta \to 0$, such that 
	\begin{equation}
	\|u^{\vec{n}}_\delta - u \|_{L^2} +	\|p^{\vec{n}}_\delta - p \|_{L^2}  + {\|\nabla \cdot u^{\vec{n}}_\delta \|_{L^2}} \leq C \delta \|f\|_{L^2}.
	\end{equation}
\end{prpstn}

\begin{proof}
	We can see from (\ref{sfns})
	\begin{align*}
	|\widehat{u}(\xi)- \widehat{u_\delta^{\vec{n}}}(\xi)| &\leq \left(\left|\frac{1}{|\lambda_\delta^{\vec{n}}(\xi)|^2} - \frac{1}{|\xi|^2}\right|  +\left| \frac{\lambda_\delta^{\vec{n}}(\xi) \overline{{\lambda_\delta^{\vec{n}}(\xi)}}^T  }{|\lambda_\delta^{\vec{n}}(\xi)|^4} -\frac{i \xi (-i \xi)^T}{|\xi|^4}\right|\right) |\widehat{f}(\xi)| \\
	&\leq 2 \left(\frac{1}{|\lambda_\delta^{\vec{n}}(\xi)|} + \frac{1}{|\xi|}\right)  \left|\frac{\lambda_\delta^{\vec{n}}(\xi)}{|\lambda_\delta^{\vec{n}}(\xi)|^2}-\frac{i\xi}{|\xi|^2} \right| |\widehat{f}(\xi)| \\
	& \leq C \left( \left|\frac{\lambda_\delta^{\vec{n}}(\xi)}{|\lambda_\delta^{\vec{n}}(\xi)|^2}-\frac{i\xi}{|\xi|^2} \right|\right)|\widehat{f}(\xi)|
	\end{align*}
	where the last inequality is due to Theorem \ref{t1}. Then since
	$$
	|\widehat{p}(\xi) - \widehat{p^{\vec{n}}_\delta}(\xi)| \leq \left|	\frac{-\overline{\lambda_\delta^{\vec{n}}(\xi)}}{|\lambda_\delta^{\vec{n}}(\xi)|^2} - \frac{i\xi}{|\xi|^2}\right||\hat{f}(\xi)| = \left|	\frac{{\lambda_\delta^{\vec{n}}(\xi)}}{|{\lambda_\delta^{\vec{n}}(\xi)}|^2} - \frac{i\xi}{| \xi|^2}\right||\widehat{f}(\xi)|,
	$$ it is sufficient to prove 
	$$ \left|\frac{\lambda_\delta^{\vec{n}}(\xi)}{|\lambda_\delta^{\vec{n}}(\xi)|^2}-\frac{i\xi}{|\xi|^2} \right| \leq C \delta \qquad (\dagger)$$ for some $C$ independent of $\delta$ and $\xi$. To this end let us explicitly write $\lambda_\delta^{\vec{n}}(\xi) = a^{\vec{n}}_\delta(\xi) + i b_\delta(\xi)$ to obtain
	\begin{align*}
	\left|	\frac{{\lambda_\delta^{\vec{n}}(\xi)}}{|{\lambda^{\vec{n}}_\delta(\xi)}|^2} - \frac{i\xi}{|\xi|^2}\right| 
	&=   \left| \frac{a^{\vec{n}}_\delta(\xi) + i b_\delta(\xi)}{|a^{\vec{n}}_\delta(\xi) + i b_\delta(\xi)|^2 } - \frac{i \xi}{|\xi|^2} \right| \\
	& \leq \underbrace{\left| \frac{ i b_\delta(\xi)}{| i b_\delta(\xi)|^2 } - \frac{i \xi}{|\xi|^2} \right|}_{I_1} + 
	\underbrace{\left|\frac{i b_\delta(\xi)}{|a^{\vec{n}}_\delta(\xi) + i b_\delta(\xi)|^2 } -  \frac{ i b_\delta(\xi)}{| i b_\delta(\xi) |^2} \right|}_{I_2} + 
	\underbrace{\left|\frac{a^{\vec{n}}_\delta(\xi)}{|a^{\vec{n}}_\delta(\xi) + i b_\delta(\xi)|^2 }\right|}_{I_3}
	\end{align*}
	and consider cases depending on the values of $k = \delta |\xi|$ as in the proof of Theorem \ref{t1}.
	\begin{enumerate}
		\item $k <  1$.  We first consider $I_1$.
		Using Lemma \ref{l2} we have
		$$
		I_1 \leq \left|\frac{1}{\Lambda_\delta(\xi)}-\frac{1}{|\xi|}\right|  = \delta \left|\frac{1}{\displaystyle 4 \int_{r=0}^{1}\int_{\theta=0}^{\frac{\pi}{2}} w(r)r\cos(\theta)\sin(k r \cos(\theta))dr d\theta} - \frac{1}{k}\right|\,.
		$$
		Since
		$$
		x  - \frac{x^3}{6}\leq \sin(x) \leq x  \text{ for } 0 \leq x \leq 1,
		$$
		we obtain 
		$$
		2\pi k - \frac{\pi k^3}{4} \leq 4 \int_{r=0}^{1}\int_{\theta=0}^{\frac{\pi}{2}} w(r)r\cos(\theta)\sin(k r \cos(\theta))dr d\theta \leq 2\pi k,
		$$ which implies
		$$
		\frac{1}{\delta} \left|\frac{1}{\Lambda_\delta(\xi)}-\frac{1}{|\xi|}\right| \leq \frac{1}{2\pi}\left(\frac{1}{k-\frac{k^3}{8}} - \frac{1}{k}\right) = \frac{1}{2\pi} \frac{k}{8-k^2} \leq \frac{k}{2\pi}.
		$$
		Hence $I_1 \leq \frac{\delta k}{2\pi} \leq \frac{\delta}{2\pi}$.
		\bigskip
		
		Next, for $I_2$, since $\cos(x) - 1 \leq \frac{x^2}{2}$, it is clear
		\begin{equation*} 
		|a^{\vec{n}}_\delta(\xi)| = \left|2 \int_{\mathcal{H}_{\vec{n}}} \frac{w_\delta(|s|)s}{|s|}(\cos(\xi\cdot s)-1) ds \right| \leq C \delta |\xi|^2.
		\end{equation*}
		where we have used the moment condition (\ref{mc}).
		
		On the other hand we can see from the proof of Theorem \ref{t1} that
		$$
		|b_\delta(\xi)| \geq C |\xi|\,.
		$$
		
		Hence it follows that 
		$$
		I_2  \leq \frac{|a^{\vec{n}}_\delta(\xi)|^2}{(|a^{\vec{n}}_\delta(\xi)|^2+|b_\delta(\xi)|^2)|b_\delta(\xi)|} \leq  \frac{C \delta^2 |\xi|^4 }{ |\xi|^3} \leq C \delta\,.
		$$
		
		As for  $I_3$, similar calculation as in the case of $I_2$ shows $I_3 \leq C \delta$.

		\item $k \geq 1$.
		We observe from the proof of Theorem \ref{t1} that $|{\lambda_\delta^{\vec{n}}(\xi)}|\geq \frac{C}{\delta}$, hence
		$$
		\frac{1}{\delta}\left|	\frac{{\lambda_\delta^{\vec{n}}(\xi)}}{|{\lambda_\delta^{\vec{n}}(\xi)}|^2} - \frac{\xi}{|\xi|^2}\right| \leq \frac{1}{\delta|\lambda_\delta^{\vec{n}}(\xi)|} + \frac{1}{\delta|\xi|} \leq {C}.$$
	\end{enumerate}
	{Lastly the local divergence of $u_\delta$ can be estimated as 
		\begin{align*}
		|\widehat{\nabla \cdot u^{\vec{n}}_\delta}(\xi)| &\leq \frac{|\xi|}{|\lambda_\delta^{\vec{n}}(\xi)|} \left(\frac{|\xi|}{|\lambda_\delta^{\vec{n}}(\xi)|}\left|	\frac{{\lambda_\delta^{\vec{n}}(\xi)}}{|{\lambda_\delta^{\vec{n}}(\xi)}|^2} - \frac{i\xi}{| \xi|^2}\right| +  \left|\frac{\lambda_\delta^{\vec{n}}(\xi)}{|\lambda_\delta^{\vec{n}}(\xi)|^2}-\frac{i\xi}{|\lambda_\delta^{\vec{n}}(\xi)|^2}\right| \right) |\widehat{f}(\xi)| \\
		& \leq \frac{|\xi|}{|\lambda_\delta^{\vec{n}}(\xi)|} \left(\frac{2|\xi|}{|\lambda_\delta^{\vec{n}}(\xi)|} + 1 \right)\left|	\frac{{\lambda_\delta^{\vec{n}}(\xi)}}{|{\lambda_\delta^{\vec{n}}(\xi)}|^2} - \frac{i\xi}{| \xi|^2}\right|  |\widehat{f}(\xi)|
		\leq C\delta |\widehat{f}(\xi)|
		\end{align*} where the last inequality is due to the estimate $(\dagger)$.}
\end{proof}

The study on the nonlocal Stokes equation can be extended to time-dependent case. Let us consider
\begin{eqnarray} \label{uns}
(u^{\vec{n}}_\delta)_t -\mathcal{L}^{\vec{n}}_\delta u^{\vec{n}}_\delta + \mathcal{G}_\delta^{\vec{n}} p^{\vec{n}}_\delta &=& f \nonumber \text{ in } (0,T) \times \Omega  \\
-\mathcal{D}^{\vec{n}}_\delta u^{\vec{n}}_\delta &=& 0 \text{ in } (0,T) \times \Omega  \nonumber \\
u^{\vec{n}}_\delta \vert_{t=0} &=& u_0 \text{ in } \Omega
\end{eqnarray}
along with its counterpart local equation
\begin{eqnarray} \label{uls}
u_t - \Delta u + \nabla p &=& f \nonumber \text{ in } (0,T) \times \Omega   \\
-\nabla \cdot u &=& 0 \text{ in } (0,T) \times \Omega \nonumber \\
u \vert_{t = 0} &=& u_0 \text{ in } \Omega
\end{eqnarray}
where all the local and nonlocal field variables as well as the data are assumed to be periodic on $\Omega$ with zero means. We then have the analogous results as in the steady case.
\begin{prpstn} \label{tns}
	Assume $f \in L^2(0,T: (S_\delta^{\vec{n}}(\Omega))^{\star})$ and $u_0 \in L^2(\Omega)$ with $-\mathcal{D}^{\vec{n}}_\delta u_0 = 0$ in $\Omega$ . Then the nonlocal Stokes equation (\ref{uns}) has a unique solution $(u^{\vec{n}}_\delta,p^{\vec{n}}_\delta)$ where $u^{\vec{n}}_\delta \in L^2(0,T: S_\delta(\Omega)) \cap C(0,T:L^2(\Omega))$, $(u^{\vec{n}}_\delta)_t \in L^2(0,T: (S_\delta^{\vec{n}}(\Omega))^{\star})$ and $p^{\vec{n}}_\delta \in L^2(0,T: L^2(\Omega))$. 
\end{prpstn}
\begin{proof}
	Let us write $P^{\vec{n}}_\delta$ to denote the nonlocal Leray operator which in the Fourier space is given by 
	$$\widehat{P^{\vec{n}}_\delta}(\xi) :=  I - \frac{\lambda^{\vec{n}}_\delta(\xi)\overline{{\lambda^{\vec{n}}_\delta(\xi)}}^T  }{|\lambda^{\vec{n}}_\delta(\xi)|^2}.$$
	{One can check that $P^{\vec{n}}_\delta u = u$ for $\mathcal{D}_\delta^{\vec{n}} u = 0$, $P^{\vec{n}}_\delta$ commutes with $\mathcal{L}_\delta^{\vec{n}}$, and $P^{\vec{n}}_\delta \circ \mathcal{G}_\delta^{\vec{n}} = 0$,  hence the nonlocal system (\ref{uns}) is equivalent to}
	\begin{eqnarray} 
	(u_\delta^{\vec{n}})_t -\mathcal{L}^{\vec{n}}_\delta u^{\vec{n}}_\delta &=& P^{\vec{n}}_\delta f \nonumber \text{ in } (0,T) \times \Omega  \\
	\mathcal{G}^{\vec{n}}_\delta p^{\vec{n}}_\delta &=& f-P^{\vec{n}}_\delta f \text{ in } (0,T) \times \Omega  \nonumber \\
	u^{\vec{n}}_\delta \vert_{t=0} &=& u_0 \text{ in } \Omega \nonumber
	\end{eqnarray}
	of which the unique solutions are given by Duhamel's principle
	$$
	\widehat{u^{\vec{n}}_\delta}(\xi,t) = \widehat{u^{\vec{n}}_0}(\xi) \exp(-|\lambda^{\vec{n}}_\delta(\xi)|^2t) + \int_{0}^{t} \exp(-|\lambda^{\vec{n}}_\delta(\xi)|^2(t-s))\widehat{P^{\vec{n}}_\delta}(\xi)\widehat{f}(\xi,s) ds 	
	$$
	and 
	$$
	\widehat{p^{\vec{n}}_\delta}(\xi,t) = \frac{\overline{\lambda^{\vec{n}}_\delta(\xi)}^T}{|\lambda^{\vec{n}}_\delta(\xi)|^2}{(I - \widehat{P^{\vec{n}}_\delta}(\xi))\widehat{f}(\xi,t)}. 	
	$$
	{We may then apply the standard energy arguments to show $u^{\vec{n}}_\delta,(u^{\vec{n}}_\delta)_t$ and $p_\delta^{\vec{n}}$ belong to the appropriate spaces}. In order to show the continuity of $u_\delta$, we first deduce from Theorem \ref{t1} that $S_\delta(\Omega) \subset L^2(\Omega) \subset S^{\star}_\delta(\Omega)$ is a Hilbert triple and then apply the classical interpolation result \cite{temam2001}. 
\end{proof}

To conclude our discussion of the nonlocal Stokes equation we prove that the nonlocal solutions of (\ref{uns}) converge to the corresponding local ones as the nonlocal parameter $\delta$ vanishes. {Given a locally divergence free initial velocity $u_0$, however, we need to exercise care in prescribing the initial velocity for the nonlocal Stokes equations }{since Lemma \ref{l2} shows that $u_0$ is in general not nonlocally divergence free.} 
\begin{prpstn} \label{rluns}
	Suppose $f \in L^2(0,T: L^2(\Omega))$ and $u_0 \in L^2(\Omega)$ {with $-\nabla \cdot u_0 = 0$ in $\Omega$. Assume $u_{0,\delta}^{\vec{n}} \to u_0$ in $L^2(\Omega)$ as $\delta \to 0$ with $-\mathcal{D}_\delta^{\vec{n}} u_{0,\delta}^{\vec{n}} = 0$. Then the unique solution $(u^{\vec{n}}_\delta,p^{\vec{n}}_\delta)$ of (\ref{uns}) where $u_0$ is replaced by $u_{0,\delta}^{\vec{n}}$} converges to the unique solution $(u,p)$ of (\ref{uls}) in $L^2(0,T:L^2(\Omega))$ as $\delta \to 0$.
\end{prpstn}  
\begin{proof}
	Let us define the local Leray projector $P_0$ in the Fourier space
	$$
	\widehat{P_0}(\xi) = I - \frac{\xi \xi^T}{|\xi|^2}.
	$$ 
	The local solutions $u$ and $p$ are then given by
	$$
	\widehat{u}(\xi,t) = \widehat{u_0}(\xi) \exp(-|\xi|^2t) + \int_{0}^{t} \exp(-|\xi|^2(t-s))\widehat{P_0}(\xi) \widehat{f}(\xi,s) ds
	$$
	and
	$$
	\widehat{p}(\xi,t) = \frac{-i{\xi}^T}{|\xi|^2}{(I - \widehat{P_0}(\xi))\widehat{f}(\xi,t)}.
	$$
	We first consider 
	{
		\begin{align*}
		\widehat{u_\delta^{\vec{n}}}(\xi,t)-\widehat{u}(\xi,t) &=  \widehat{u_{0,\delta}^{\vec{n}}}(\xi)\exp(-|\lambda^{\vec{n}}_\delta(\xi)|^2t)-\widehat{u_0}(\xi)\exp(-|\xi|^2t)  \\ & + \int_{0}^{t} {\left( \exp(-|\lambda^{\vec{n}}_\delta(\xi)|^2(t-s))\widehat{P}^{\vec{n}}_\delta(\xi)-\exp(-|\xi|^2(t-s))P_0(\xi) \right) \widehat{f}(\xi,s)} ds
		\end{align*}
		where $\widehat{P^{\vec{n}}_\delta}(\xi)$ is the nonlocal Leray operator used in the proof of Proposition \ref{tns}, hence 
		\begin{align*}
		\int_0^{T}\|u_\delta-u\|_2^2 dt &\leq C_T {\int_0^{T}} \left(\sum_{\xi \in \mathbb{Z}^2, \xi \neq 0}   \left(|\widehat{u_{0}}(\xi)|^2 + |\widehat{u^{\vec{n}}_{0,\delta}}(\xi)|^2\right) \left| \exp(-|\lambda^{\vec{n}}_\delta(\xi)|^2t) - \exp(-|\xi|^2t)\right|^2 \right. + \nonumber \\ & \left. \int_{0}^{t} \left|\left(\exp(-|\lambda^{\vec{n}}_\delta(\xi)|^2(t-s))\widehat{P}^{\vec{n}}_\delta(\xi)-\exp(-|\xi|^2(t-s))\widehat{P}_0(\xi)\right) \widehat{f}(\xi,s)\right|^2 ds \right) dt
		\end{align*}}
	where $C_T$ is a constant depending only on $T$.
	We observe that the integrand in the parentheses is bounded, uniformly in $\delta$ and $t$, by $C(\|u_0\|_2^2 + \|f\|^2_{L^2(0,T:L^2(\Omega))})$. {One can easily verify 
		$$
		|\lambda^{\vec{n}}_\delta(\xi) - i\xi| \to 0, \text{ hence }
		\left|\widehat{P}^{\vec{n}}_\delta(\xi)-\widehat{P}_0(\xi)\right| \leq \left|\widehat{P}^{\vec{n}}_\delta(\xi)-\widehat{P}_0(\xi)\right|_F \to 0 
		\text{ as } \delta \to 0 \text{ for each } \xi \in \mathbb{Z}^2, \xi \neq 0$$
		so that}
	it follows the dominated convergence theorem yields  $u_\delta \to u $ in $L^2(0,T:L^2(\Omega))$ as $\delta \to 0$. 
	
	We apply the similar argument to the expression
	$$
	\int_0^{T} \|{p^{\vec{n}}_\delta}-p\|_2^2 dt = \int_0^T \sum_{\xi \in \mathbb{Z}^2, \xi \neq 0} \left( \frac{\overline{\lambda^{\vec{n}}_\delta(\xi)}^T}{|\lambda^{\vec{n}}_\delta(\xi)|^2}{(I - \widehat{P_\delta}(\xi))\widehat{f}(\xi,t)} + \frac{i{\xi}^T}{|\xi|^2}{(I - \widehat{P_0}(\xi))\widehat{f}(\xi,t)}\right)^2 dt
	$$
	to conclude  $p^{\vec{n}}_\delta \to p$ in $L^2(0,T:L^2(\Omega))$ as $\delta \to 0$.
	
	{Lastly we remark that nonlocally divergence free initial velocity $u_{0,\delta}$ as in the assumption of the theorem can be explicitly constructed by taking $\widehat{u_{0,\delta}^{\vec{n}}}(\xi) = \widehat{P^{\vec{n}}_\delta}(\xi) \widehat{u_0}(\xi)$.}
\end{proof}

We note that one can also  get the order of convergence as in  the time-independent case. The details are omitted.

\subsection{Nonlocal Helmholtz decomposition}
{The nonlocal Leray operators introduced in the proof of Proposition \ref{tns} clearly imply a nonlocal version of the classical Helmholtz decomposition theorem which warrants a more detailed discussion.} 

We begin with the following two dimensional result.
\begin{thrm}		\label{2dnhd}
	If $u:\Omega \to \mathbb{R}^2$ is periodic with zero mean, then there exist unique periodic, zero mean scalar potentials $p_\delta^{\vec{n}}$ and $q_\delta^{\vec{n}}$ on $\Omega$ such that
	$$
	u(x) = \mathcal{G}_\delta^{\vec{n}} p_{\delta}^{\vec{n}}(x) + 		\begin{bmatrix}
	0 & -1 \\
	1 & 0 
	\end{bmatrix}\mathcal{G}_\delta^{-\vec{n}} q_\delta^{\vec{n}}(x)
	$$ 
	with 	$ \mathcal{D}_\delta^{\vec{n}}\left( \begin{bmatrix}
	0 & -1 \\
	1 & 0 
	\end{bmatrix}\mathcal{G}_\delta^{-\vec{n}} q_\delta^{\vec{n}}\right)(x) = 0.$ {In addition we have  the estimate $$
		\|p^{\vec{n}}_\delta\|_{S^{\vec{n}}_\delta(\Omega)} + \|q^{\vec{n}}_\delta\|_{S^{\vec{n}}_\delta(\Omega)} \leq {C} \|u\|_{L^2(\Omega)}
		$$ for some constant {$C$} independent of $\vec{n}$ and also of $\delta$ as $\delta \to 0$. Here $S_\delta^{\vec{n}}(\Omega)$ is the energy space as in Lemma \ref{wns}. }
\end{thrm}
\begin{proof}
	In the Fourier space the unique solutions are given by
	$$\widehat{p_\delta^{\vec{n}}}(\xi) = {-}\frac{{\lambda_\delta^{-\vec{n}}(\xi)}^T \widehat{u}(\xi)}{|\lambda_\delta^{\vec{n}}(\xi)|^2} \quad \text{ and } \quad \widehat{q_\delta^{\vec{n}}}(\xi) = \frac{{\lambda_\delta^{\vec{n}}(\xi)}^T}{|\lambda_\delta^{-\vec{n}}(\xi)|^2}\begin{bmatrix}
	0 & 1 \\
	-1 & 0 
	\end{bmatrix} {
		\left(I {+} \frac{{\lambda_\delta^{\vec{n}}(\xi)}{\lambda_\delta^{-\vec{n}}(\xi)}^T}{|\lambda^{\vec{n}}_\delta(\xi)|^2}\right)
	} \widehat{u}(\xi)$$ and the rest of the claim is due to ${{\lambda_\delta^{-\vec{n}}(\xi)}} = -\overline{\lambda_\delta^{\vec{n}}(\xi)}$ and Theorem \ref{t1}.
\end{proof}
In three dimensions we introduce nonlocal curl operators
$$
\mathcal{C}_\delta^{\vec{n}}v(x)  = 2\int_{\mathbb{R}^3} {\chi}_{\vec{n}}(y-x)  \underline{\omega}_\delta(y-x) (y-x) \times \frac{v(y)-v(x)}{|y-x|} dy,
$$
based on which we deduce the following.
\begin{thrm} \label{3dnhd}
	If $u: \widetilde{\Omega}: (-\pi,\pi)^3 \to \mathbb{R}^3$ is periodic with zero mean, then there exist unique periodic, zero mean scalar and vector potentials $p_\delta^{\vec{n}}$ and $v_\delta^{\vec{n}}$, respectively, such that
	$$
	u(x) = \mathcal{G}_\delta^{\vec{n}} p_\delta^{\vec{n}}(x) + \mathcal{C}_\delta^{-\vec{n}}v_\delta^{\vec{n}}(x)
	$$
	with the nonlocal Gauge condition $\mathcal{D}_\delta^{-\vec{n}}v_\delta^{\vec{n}}(x) = 0$. Moreover
	{ $(\mathcal{C}_\delta^{\vec{n}}\circ \mathcal{G}_\delta^{\vec{n}}) p_\delta^{\vec{n}}(x)$ and $ (\mathcal{D}_\delta^{\vec{n}} \circ \mathcal{C}_\delta^{-\vec{n}})v_\delta^{\vec{n}}(x)$ vanish along with the analogous estimate as in Theorem \ref{2dnhd} with $v_\delta^{\vec{n}}$ in place of $p_\delta^{\vec{n}}$.}
\end{thrm}
\begin{proof}
	One can verify that a nonlocal vector identity 
	$$
	(\mathcal{C}_\delta^{-\vec{n}}\circ \mathcal{C}_\delta^{\vec{n}})f(x) = (\mathcal{G}_\delta^{\vec{n}} \circ \mathcal{D}_\delta^{\vec{n}})f(x) - \mathcal{L}_\delta^{\vec{n}}f(x)
	$$
	holds for any periodic $f:(-\pi,\pi)^3 \to \mathbb{R}^3$. Solving for $f$ in $-\mathcal{L}_\delta^{\vec{n}}f = u$ then yields the Fourier representations of the unique solutions
	$$
	\widehat{p}^{\vec{n}}(\xi) =  -\frac{{\lambda_\delta^{-\vec{n}}(\xi)}^T \widehat{u}(\xi)}{|\lambda_\delta^{\vec{n}}(\xi)|^2} \quad \text{ and } \quad \widehat{v}^{\vec{n}}(\xi) = \frac{\lambda_\delta^{\vec{n}}(\xi)}{|\lambda_\delta^{\vec{n}}(\xi)|^2} \times \widehat{u}(\xi).
	$$ 
	{We refer to the proof of Theorem \ref{2dnhd} for the rest of the claim.}
\end{proof}

Let us point out that the proof of Theorem \ref{3dnhd} reveals the well-posedness of a nonlocal version of the classical first order div-curl elliptic system
as stated below.

\begin{thrm}		\label{2divcurl}
	Given periodic, zero mean data $f,g$ in $L^2(\widetilde{\Omega})$ with $\mathcal{D}_\delta^{-\vec{n}}g = 0$, 
	there exist a unique periodic, zero mean vector field 
	$u:\widetilde{\Omega} \to \mathbb{R}^{3}$ satisfying the following nonlocal div-curl system
	$$\begin{array}{ll}
	\mathcal{D}_\delta^{\vec{n}}u = f,  & \quad \mbox{in}\; \widetilde{\Omega},\\
	\mathcal{C}_\delta^{\vec{n}}u = g, &\quad \mbox{in}\; \widetilde{\Omega},
	\end{array}$$
	and the estimate
	\begin{equation} \label{2dfried}		
	\int_{\widetilde{\Omega}} \left(|u(x)|^2 + |\mathcal{G}_\delta^{\vec{n}}u(x)|^2 \right)dx \leq C \int_{\widetilde{\Omega}} \left(|\mathcal{D}_\delta^{\vec{n}} u(x)|^2 + |\mathcal{C}_\delta^{\vec{n}} u(x)|^2 \right)dx 
	\end{equation} for some positive constant $C$ independent of $u$, $\vec{n}$ and $\delta$ (as $\delta \to 0$).
\end{thrm}
The estimate \eqref{2dfried}  is  a nonlocal version of the second Friedrichs inequality \cite{kn84,s83}, which can be easily shown with the help of Theorem \ref{t1}. We omit the details.

We note that classical Helmholtz decomposition have wide applications of mechanics and electromagnetics.
For nonlocal and fractional versions, one can also check \cite{Tar08bk}. The further study of nonlocal div-curl systems, including nonlocal versions of the div-curl lemma \cite{cdm11div,murat78}, may also be of interests and will be left for future works.

\subsection{Nonlocal correspondence models of isotropic linear elasticity}
We study the elastic potential energy given by  
\begin{equation} \label{nep}
{\mathcal{E}_\delta^{\vec{n}}(\mathbf{u}) = \frac{1}{2}\lambda \|{\mathcal{D}_\delta^{\vec{n}}\mathbf{u(x)}}\|^2_2 + \mu \|e^{\vec{n}}_\delta(\mathbf{u(x)})\|^2_2}
\end{equation} 	 
where $\lambda,\mu$ are Lam\'e coefficients and $e^{\vec{n}}_\delta(\mathbf{u})$ is the nonlocal strain tensor 
$$ e^{\vec{n}}_\delta(\mathbf{u}) = \frac{\mathcal{G}^{\vec{n}}_\delta \mathbf{u}+ (\mathcal{G}^{\vec{n}}_\delta \mathbf{u})^T }{2}$$ for a displacement field $\mathbf{u}:\Omega \to \mathbb{R}^2$. 
This can be viewed as a so-called nonlocal correspondence model where the local stress tensors in the  local energy density are replaced by the nonlocal counterparts \cite{dt17,s17}. We assume $\mu>0$ and ${\lambda} + {2\mu}>0$ for which it is well known that the corresponding local elastic energy is variationally stable over the energy space 
$${V}_0 = \{\mathbf{u}:\Omega \to \mathbb{R}^2  \vert  \text{ each } u_i \text{ periodic with zero mean } i = 1, 2 \text{ and }  \sum_{i=1}^{2} \|u_i\|^2_{2} + \|\nabla u_i\|^2_{2} < \infty \}.$$ 
The energy $\mathcal{E}_\delta^{\vec{n}}$ can be seen as a linear elastic potential energy of isotropic materials in the peridynamics correspondence theory \cite{s00,s17}. As in the one dimensional case with radially symmetric interaction kernels \cite{dt17}, let us first define the nonlocal energy space $V_\delta^{\vec{n}}$ to be the closure of $C^\infty$, zero mean, periodic $\mathbb{R}^2$-valued functions on $\Omega$ with respect to the norm 
$$\|\mathbf{u}\|_{V_\delta^{\vec{n}}} = (\|\mathbf{u}\|^2_2 + \mathcal{E}_\delta^{\vec{n}}(\mathbf{u}))^{1/2}.$$
A precise statement of the variational stability of $\mathcal{E}_\delta^{\vec{n}}$ is then given by
\begin{equation} \label{vs}
\mathcal{E}^{\vec{n}}_\delta(\mathbf{u}) \ge C \|\mathbf{u}\|^2_{V_\delta} \qquad \forall \mathbf{u} \in V_\delta^{\vec{n}}
\end{equation}
for some constant $C>0$ independent of $\vec{n}$ and $\delta$, as $\delta \to 0$. {In order to establish the stability let us introduce the nonlocal Navi\'er operator $P_\delta^{\vec{n}}$ defined as
	\begin{eqnarray}
	P^{\vec{n}}_\delta(\mathbf{u}) = -\mu \mathcal{L}_\delta^{\vec{n}} \mathbf{u} - (\lambda+\mu) \mathcal{G}_\delta^{\vec{n}} (\mathcal{D}_\delta^{\vec{n}} \mathbf{u})
	\end{eqnarray}
	so that for $\mathbf{u} = (u_1,u_2)\in V_\delta^{\vec{n}}$
	\begin{align*}
	(P^{\vec{n}}_\delta(\mathbf{u}),\mathbf{u})_2 
	& = \lambda \int_{\Omega} |\mathcal{D}_\delta^{\vec{n}} \mathbf{u}(\mathbf{x})|^2 d\mathbf{x}  + \mu \int_{\Omega} \left( |\mathcal{D}_\delta^{\vec{n}} \mathbf{u}(\mathbf{x})|^2  + \sum_{i=1,2} |\mathcal{G}_\delta^{\vec{n}} u_i(\mathbf{x})|^2  \right) d\mathbf{x} = 2 \mathcal{E}_\delta^{\vec{n}}(\mathbf{u})
	\end{align*} due to \eqref{anib}. We observe from Lemma \ref{l1} $\mathcal{D}_\delta^{\vec{n}}\mathbf{u} =  \text{Tr}(e^{-\vec{n}}_\delta(\mathbf{u}))$ is not in general equal to $\text{Tr}(e^{\vec{n}}_\delta(\mathbf{u}))$ where Tr denotes the trace operator, which is different from the local case. } 

Similar to the studies on the classical Korn's inequality and the nonlocal versions in \cite{k1909,me12ccm,nh17}, we
have a nonlocal version for the nonlocal Navi\'er system as follows.
\begin{lmm}{(Nonlocal Korn's inequality).} \label{nck}
	There exists a constant $C>0$ independent of $\vec{n}$ and $\delta$ such that
	$$
	\mathcal{E}^{\vec{n}}_\delta(\mathbf{u}) \geq 	 C \|\mathcal{G}_{\delta}^{\vec{n}} \mathbf{u}(\mathbf{x})\|^{{2}}_2  \qquad \forall \mathbf{u} \in V_\delta^{\vec{n}} .
	$$ 	
\end{lmm}
\begin{proof}
	We can see
	\begin{align*}
	\mathcal{E}^{\vec{n}}_\delta(\mathbf{u}) =  \sum_{\xi \in \mathbb{Z}^2, \xi \neq (0,0)} \widehat{\mathbf{u}}(\xi) \cdot P_\delta^{\vec{n}}(\xi) \widehat{\mathbf{u}}(\xi) &= \sum_{\xi \in \mathbb{Z}^2, \xi \neq (0,0)} \mu |\lambda_\delta^{\vec{n}}(\xi)|^2|\widehat{\mathbf{u}}(\xi)|^2 + (\mu+\lambda) |\lambda_\delta^{\vec{n}}(\xi)\cdot \widehat{\mathbf{u}}(\xi) |^2 \\
	& \geq \min(\mu,\lambda + 2\mu) \sum_{\xi \in \mathbb{Z}^2, \xi \neq (0,0)} |\lambda_\delta^{\vec{n}}(\xi)|^2|\widehat{\mathbf{u}}(\xi)|^2  
	\end{align*} {which proves the claim.} 
\end{proof} 
\noindent We point out that by setting $\lambda=0$, one can recover the periodic versions of nonlocal  Korn's inequality studied in  \cite{me12ccm,nh17}.

Now using  Theorem \ref{t1}, which serves like a nonlocal Poincare inequality \cite{du19,md16},  we readily have the following coercivity result and the well-posedness.
\begin{lmm} \label{nc}
	There exists a constant $C>0$ independent of $\vec{n}$ and $\delta$ as $\delta \to 0$ such that
	$$
	\mathcal{E}^{\vec{n}}_\delta(\mathbf{u}) \geq
	C \| \mathbf{u} \|_2^2
	\qquad \forall \mathbf{u} \in V_\delta^{\vec{n}}.
	$$ 	
	Moreover,  the problem
	\begin{equation} \label{nle}
	P_\delta^{\vec{n}}(\mathbf{u}) = \mathbf{f}
	\end{equation}
	is well-posed
	{over $V_\delta^{\vec{n}}$} where $\mathbf{f}\in L^2$. 
\end{lmm}
We note that the same remains valid for $\mathbf{f}$ 
belonging to the dual space $(V_\delta^{\vec{n}})^*$. Indeed, given the equivalence of $\mathcal{E}_\delta^{\vec{n}}(\cdot)$ with $\|\cdot \|_{V_\delta^{\vec{n}}}$ due to Lemma \ref{nc},
an explicit characterization of the dual space with the $L^2$ duality pairing {can also be obtained} as done in \cite{zd10}. 
{Not only is the nonlocal solution $\mathbf{u}_\delta^{\vec{n}} \in V_\delta^{\vec{n}}$ to (\ref{nle}) important in its own right as a unique minimizer  of $
	\mathcal{E}_\delta^{\vec{n}}(\mathbf{u})-(\mathbf{f},\mathbf{u})_2
	$ but it can also be shown to recover the corresponding local solution to the local Navi\'er equation, when the latter is well-posed. To this end let us first present the following embedding result.}

\begin{lmm} \label{ue}
	There exists a constant $C$ independent of $\delta$ and $\vec{n}$ such that
	$$
	\|\mathbf{u}\|_{V_\delta^{\vec{n}}} \leq C \|\mathbf{u}\|_{V_0} \quad \forall \mathbf{u} \in V_0.
	$$
\end{lmm}

\begin{proof}
	{As can be observed in the proof of Lemma \ref{nck},  we have $$\widehat{\mathbf{u}}(\xi) \cdot P^{\vec{n}}_\delta(\xi) \widehat{\mathbf{u}}(\xi) \leq \max(\mu,\lambda + 2\mu) |\lambda_\delta^{\vec{n}}(\xi)|^2|\widehat{\mathbf{u}}(\xi)|^2 $$ to which applying Theorem \ref{t1} proves the claim.} 
\end{proof}

{It is in this nonlocal norm $\|\cdot\|_{{V}_{\delta}^{\vec{n}}}$, hence in $\|\cdot\|_{L^2}$, that we have the convergence of the nonlocal solution $\mathbf{u}_\delta$ to its local counterpart. }

\begin{prpstn}
	Let $\mathbf{u}_\delta^{\vec{n}}$ denote the solution of the nonlocal Navi\'er equation (\ref{nle}). Then there exists a constant $C$ independent of $\delta$ and $\vec{n}$ as $\delta \to 0$ such that
	$$\|\mathbf{u}_\delta^{\vec{n}} - \mathbf{u}\|_{{{V^{\vec{n}}_{\delta}}}} \leq C \delta \|\mathbf{f}\|_2$$
	where $\mathbf{u}$ is the solution to the local Navi\'er equation
	$$
	P_0(\mathbf{u}):= -\mu \Delta \mathbf{u} - (\lambda+\mu) \nabla \nabla \cdot \mathbf{u} = \mathbf{f} \quad \text{ in } \Omega
	$$
	subject to the periodic boundary condition.
\end{prpstn}

\begin{proof}
	{Let us first write
		\begin{align*}
		\|\mathbf{u}_\delta^{\vec{n}} - \mathbf{u}\|^2_{{{V}_\delta^{\vec{n}}}} &= \sum_{\xi \in \mathbb{Z}^2, \xi \neq (0,0)} \widehat{\mathbf{f}}(\xi) \cdot ((\widehat{P}_\delta^{\vec{n}}(\xi))^{-1}- (\widehat{P}_0(\xi))^{-1})(I-\widehat{P}_\delta^{\vec{n}}(\xi)(\widehat{P}_0(\xi))^{-1})\widehat{\mathbf{f}}(\xi) \\
		& \leq C \sum_{\xi \in \mathbb{Z}^2, \xi \neq (0,0)}  |(\widehat{P}_\delta^{\vec{n}}(\xi))^{-1}- (\widehat{P}_0(\xi))^{-1}||\widehat{\mathbf{f}}(\xi)|^2
		\end{align*}
		where the inequality is due to Lemma \ref{ue}. But one can apply Theorem \ref{t1} and 
		the following estimate shown in the proof of Proposition \ref{p1}
		$$ \left|\frac{1}{|\lambda_\delta(\xi)|}-\frac{1}{|\xi|} \right| \leq C \delta $$ to verify that
		$$|(\widehat{P}_\delta^{\vec{n}}(\xi))^{-1}- (\widehat{P}_0(\xi))^{-1}|_F \leq C \delta $$ where $|\cdot|_F$ denotes the Frobenius norm, from which the result follows.}
\end{proof}

As in the previous section on the nonlocal Stokes equation, we now consider the time dependent nonlocal Navi\'er equation
\begin{align}  \label{tdnne}
(\mathbf{u}^{\vec{n}}_\delta)_{tt} + P_\delta^{\vec{n}}\mathbf{u}^{\vec{n}}_\delta &= \mathbf{f} \nonumber \text{ in } (0,T) \times \Omega\\
\mathbf{u}^{\vec{n}}_\delta \vert_{t=0} &= \mathbf{g} \nonumber \text{ in } \Omega \\
(\mathbf{u}^{\vec{n}}_\delta)_t \vert_{t=0} &= \mathbf{h} \text{ in } \Omega 
\end{align}
in juxtaposition with the local Navi\'er equation
\begin{align}  \label{tdlne}
\mathbf{u}_{tt} + P_0\mathbf{u} &= \mathbf{f} \nonumber \text{ in } (0,T) \times \Omega\\
\mathbf{u} \vert_{t=0} &= \mathbf{g} \nonumber \text{ in } \Omega \\
\mathbf{u}_t \vert_{t=0} &= \mathbf{h} \text{ in } \Omega.
\end{align}

Since the Hermitian matrix $P_\delta^{\vec{n}}$ is positive definite, we can apply the similar argument as in \cite{dz11} to establish the following well-posedness result for which we omit the details.

\begin{prpstn} \label{wpusneq}
	Suppose $\mathbf{g}\in V_\delta^{\vec{n}}, \mathbf{h} \in L^2(\Omega)$ and $\mathbf{f} \in L^2(0,T:L^2(\Omega))$ where {$\mathbf{h}$ and $\mathbf{g}$} are periodic on $\Omega$ with zero means. Then there exists a unique $\mathbf{u}^{\vec{n}}_\delta$ to (\ref{tdnne}) such that 
	$$\mathbf{u}^{\vec{n}}_\delta \in C(0,T:V^{\vec{n}}_\delta), \quad (\mathbf{u}^{\vec{n}}_\delta)_t \in L^2(0,T:L^2(\Omega)).$$
\end{prpstn}

For completeness we consider convergence of the {time dependent} nonlocal solution $\mathbf{u}^{\vec{n}}_\delta$ to {the corresponding local solution} as $\delta \to 0$. {As in the steady case this can be readily established using the explicit Fourier representations of both nonlocal and local solutions.}

\begin{prpstn}
	Suppose $\mathbf{g}\in V_0(\Omega), \mathbf{h} \in L^2(\Omega)$ and $\mathbf{f} \in L^2(0,T:L^2(\Omega))$ where {$\mathbf{h}$ and $\mathbf{g}$} are periodic on $\Omega$ with zero means. Let $\mathbf{u}^{\vec{n}}_\delta$ and $\mathbf{u}$ denote the solutions of nonlocal and local Navi\'er equations, respectively, with the same initial displacement field {$\mathbf{g}$} and velocity field $\mathbf{h}$. Then we have $\mathbf{u}^{\vec{n}}_\delta \to \mathbf{u}$ in {$L^2(0,T:V_\delta^{\vec{n}}(\Omega)) \bigcap H^1(0,T:L^2(\Omega))$}.
\end{prpstn}

\begin{proof}
	{We use the explicit Fourier representation of the solutions as given in the proof of Theorem 2.25 in \cite{dz11} to write
		\begin{align*}
		\widehat{\mathbf{u}^{\vec{n}}_\delta}(t,\xi) - \widehat{\mathbf{u}}(t,\xi) &= \underbrace{\left(\cos\left(\sqrt{\widehat{P}_\delta^{\vec{n}}(\xi)}t\right) - \cos\left(\sqrt{\widehat{P}_0(\xi)}t\right)\right)\widehat{\mathbf{g}}(\xi)}_{\widehat{\mathbf{u}}_{c1}} + \underbrace{\left(\frac{\sin\left(\sqrt{\widehat{P}_\delta^{\vec{n}}(\xi)}t\right)}{\sqrt{\widehat{P}_\delta^{\vec{n}}(\xi)}} - \frac{\sin\left(\sqrt{\widehat{P}_0(\xi)}t\right)}{\sqrt{\widehat{P}_0(\xi)}}\right)\widehat{\mathbf{h}}(\xi)}_{\widehat{\mathbf{u}}_{c2}} \\
		& + \underbrace{\int_0^{t} \left(\frac{\sin\left(\sqrt{\widehat{P}_\delta^{\vec{n}}(\xi)}s\right)}{\sqrt{\widehat{P}_\delta^{\vec{n}}(\xi)}} - \frac{\sin\left(\sqrt{\widehat{P}_0(\xi)}s\right)}{\sqrt{\widehat{P}_0(\xi)}}\right) \widehat{\mathbf{f}}(t-s,\xi)ds}_{\widehat{\mathbf{u}}_{c3}}.
		\end{align*}
		We cannot argue as in \cite{dz11} since the operators $P_\delta^{\vec{n}}$ and $P_0$ do not commute, hence not simultaneously diagonalizable. Instead we argue directly as follows. We have
		\begin{align*}
		&\left(\cos\left(\sqrt{\widehat{P}_\delta^{\vec{n}}(\xi)}t\right) - \cos\left(\sqrt{\widehat{P}_0(\xi)}t\right)\right)\widehat{\mathbf{g}}(\xi) \cdot \widehat{P}_\delta^{\vec{n}}(\xi)\left(\cos\left(\sqrt{\widehat{P}_\delta^{\vec{n}}(\xi)}t\right) - \cos\left(\sqrt{\widehat{P}_0(\xi)}t\right)\right)\widehat{\mathbf{g}}(\xi) \\
		& \leq C \left(\cos\left(\sqrt{\widehat{P}_\delta^{\vec{n}}(\xi)}t\right) - \cos\left(\sqrt{\widehat{P}_0(\xi)}t\right)\right)\widehat{\mathbf{g}}(\xi) \cdot \widehat{P}_0(\xi)\left(\cos\left(\sqrt{\widehat{P}_\delta^{\vec{n}}(\xi)}t\right) - \cos\left(\sqrt{\widehat{P}_0(\xi)}t\right)\right)\widehat{\mathbf{g}}(\xi) \\
		&\leq 2\widehat{\mathbf{g}}(\xi)\cdot \widehat{P}_0(\xi)2\widehat{\mathbf{g}}(\xi)
		\end{align*}
		where the second inequality is due to Lemma \ref{ue}, hence
		$
		\|{\mathbf{u}}_{c1}\|^2_{V_\delta^{\vec{n}}} \leq 4\|\mathbf{g}\|^2_{V_0}
		$	
		uniformly in $t\in[0,T]$, which in turn implies
		$$
		\int_0^T \|{\mathbf{u}}_{c1}\|^2_{V_\delta^{\vec{n}}} dt \leq 4T\|\mathbf{g}\|^2_{V_0}.
		$$
		Let us now note that for each fixed $\xi \in \mathbb{Z}^2, \xi \neq (0,0) $ and $t \in (0,T)$ we get
		$$
		\left|\cos\left(\sqrt{\widehat{P}_\delta^{\vec{n}}(\xi)}t\right) - \cos\left(\sqrt{\widehat{P}_0(\xi)}t\right)\right| \leq \left|\cos\left(\sqrt{\widehat{P}_\delta^{\vec{n}}(\xi)}t\right) - \cos\left(\sqrt{\widehat{P}_0(\xi)}t\right)\right|_F \to 0 \quad \text{ as } \delta \to 0
		$$
		since it can be easily checked that $|\lambda^{\vec{n}}_\delta(\xi) - i\xi| \to 0$, hence $\left|\widehat{P}_\delta^{\vec{n}}(\xi) - {\widehat{P}_0(\xi)}\right|_F \to 0$, as $\delta \to 0$. The dominated convergence theorem then implies
		$$
		\int_0^T \|{\mathbf{u}}_{c1}\|^2_{V_\delta^{\vec{n}}} dt \to 0 \quad \text{ as } \delta \to 0.
		$$ Similarly the rest of the claims can be proved and we omit the details.}
\end{proof}

\section{Conclusion}
We have studied nonlocal gradient operators $\mathcal{G}_\delta^{\vec{n}}$ wherein the support of a positive kernel is prescribed to be any half-sphere parameterized by a unit vector $\vec{n}$.  This can be seen as extensions of the 
one-sided nonlocal derivative operators for scalar functions of a single variable.
It is interesting to observe that our nonlocal gradient operators with nonspherical interaction neighborhood can be effectively applied to model inherently symmetric phenomena as illustrated in our study of nonlocal Navi\'er equation of isotropic linear elasticity and nonlocal Stokes models of incompressible viscous flows. We also remark that by removing any singular growth assumption on the kernels, our nonsymmetric gradient operators are well suited to numerical quadrature based discretizations.

This work demonstrates that the symmetry of the nonlocal interaction neighborhood is not essential for nonlocal modeling
and the related mathematical theory.
While the analysis is focused on the half-sphere case, one may study further extensions that may involve only sectors of
the sphere such as those used \cite{tjd17} for nonlocal convection and in the studies of \cite{bci08,fkv15,dmt18}  on
nonlocal variational problems. The analytical results in this work are largely based on
the Fourier analysis which is limited to problems defined over periodic cells. It will be interesting to consider the analogue
on more general domains with more general boundary conditions or nonlocal constraints.
Meanwhile we point out that the well-posedness of our nonlocal Stokes and Navi\'er equations in the periodic setting is naturally linked to consideration of their Fourier spectral discretizations and related numerical issues such as the asymptotic compatibility \cite{td13} as in \cite{dtsph17}. Further numerical analysis of other discretizations are also important for applications and will be left to future works.

\begin{acknowledgement}
	\emph{Acknowledgements: }The authors would like to thank the members of CM3 group at Columbia University for discussions.
	\end{acknowledgement}
%%-----------------------------
\bibliographystyle{plain}
%\bibliography{references}

%%-----------------------------
\end{document}